\documentclass{tran-l}
\usepackage{amssymb,times}
\usepackage[numbers]{natbib}

\newtheorem*{Thm*}{Theorem}
\newtheorem{Thm}{Theorem}
\newtheorem{Cor}[Thm]{Corollary}
\newtheorem{Lemma}[Thm]{Lemma}
\newtheorem{Prop}[Thm]{Proposition}

\theoremstyle{definition}
\newtheorem{Defn}{Definition}

\newtheorem{Remark}{Remark}
\newtheorem{Ex}[Remark]{Example}

\newcommand{\mf}[1]{\mathbb{#1}}

\newcommand{\mb}[1]{\mathbf{#1}}

\DeclareMathOperator{\NC}{\mathit{NC}}

\newcommand{\abs}[1]{\left\vert#1\right\vert}

\newcommand{\set}[1]{\left\{#1\right\}}
\newcommand{\ip}[2]{\left \langle #1, #2 \right \rangle}
\newcommand{\state}[1]{\varphi \left[ #1 \right]}
\renewcommand{\phi}{\varphi}

\newcommand{\br}{\medskip\noindent}

\allowdisplaybreaks[1]

\title[Orthogonal polynomials]{Orthogonal polynomials with a resolvent-type \\
generating function}
\author[M.~Anshelevich]{Michael Anshelevich}
\thanks{This work was supported in part by an NSF grant DMS-0400860}
\address{Department of Mathematics, University of California, Riverside, CA 92521-0135}
\curraddr{Department of Mathematics, Texas A\&M University, College Station, TX 77843-3368}
\email{manshel@math.tamu.edu}
\date{October 21, 2004 and, in revised form, June 7, 2006}
\copyrightinfo{2006}
    {American Mathematical Society}

\begin{document}

\begin{abstract}
The subject of this paper are polynomials in multiple \emph{non-commuting} variables. For polynomials of this type orthogonal with respect to a state, we prove a Favard-type recursion relation. On the other hand, free Sheffer polynomials are a polynomial family in non-commuting variables with a resolvent-type generating function. Among such families, we describe the ones that are orthogonal. Their recursion relations have a more special form; the best way to describe them is in terms of the \emph{free cumulant generating function} of the state of orthogonality, which turns out to satisfy a type of second-order difference equation. If the difference equation is in fact first order, and the state is tracial, we show that the state is necessarily a rotation of a free product state. We also describe interesting examples of non-tracial infinitely divisible states with orthogonal free Sheffer polynomials.
\end{abstract}

\maketitle

\section{Introduction}
\noindent
Let $\mb{x} = (x_1, \ldots, x_n)$, $\mb{z} = (z_1, \ldots, z_n)$ be $n$-tuples of non-commuting indeterminates, such that the $x$-variables commute with the $z$-variables. Sometimes we will treat such $n$-tuples as vectors, in which case $\mb{x} \cdot \mb{z}$ denotes the scalar product.

\begin{Defn}
\label{Defn:Free-Sheffer}
Let
\[
F(\mb{z}) = 1 + \textsl{higher-order terms}
\]
be a formal power series, and $\mb{V}$ be an $n$-tuple of formal power series,
\[
V_i(\mb{z}) = z_i + \textsl{higher-order terms}.
\]
Expand
\[
F(\mb{z}) \Bigl( 1 - \mb{x} \cdot \mb{V}(\mb{z}) \Bigr)^{-1}
\]
into a power series in $z$. The coefficient of the monomial $z_{\vec{u}}$ is easily seen to be a monic polynomial $P_{\vec{u}}(\mb{x})$. We call $\set{P_{\vec{u}}}$ the (multivariate) \emph{free Sheffer polynomials}.
\end{Defn}

\br
The question to be investigated in this paper is: \emph{when are the free Sheffer polynomials orthogonal with respect to some state $\phi$?} We emphasize that $\set{x_1, \ldots, x_n}$ do not commute, and so one can talk about orthogonality of $\set{P_{\vec{u}}}$ only with respect to a positive functional on the algebra of non-commutative polynomials $\mf{R} \langle \mb{x} \rangle$, not with respect to a measure on $\mf{R}^n$.

\br
The rest of the introduction explains the motivation behind this question.

\br
Let $\set{P_n(x)}$ be a monic polynomial family (in on variable) with a generating function of the form
\begin{equation}
\label{One-variable-free-Sheffer}
\sum_{n=0}^\infty P_n(x) z^n = \frac{1}{1 - x U(z) + R(U(z))}.
\end{equation}
Here $U = z + $ \textsl{higher-order terms} and $R = z^2 + $ \textsl{higher-order terms} are formal power series. The main theorem of Section 3 of \cite{AnsMeixner} can be reformulated as follows.

\begin{Prop}
\label{Prop:One-variable-free-Meixner}
The polynomials with the generating function~\eqref{One-variable-free-Sheffer} are orthogonal if and only if the following two conditions hold:
\begin{enumerate}
\item
$U(z) = (R(z)/z)^{<-1>}$, and
\item
$R(z)/z^2 = 1 + b R(z)/z + c (R(z)/z)^2$.
\end{enumerate}
\end{Prop}

\noindent
Here $F^{<-1>}$ denotes the inverse under composition.

\br
Notice the similarity of this result to the following theorem, found in various forms by various people and going back to Meixner~\cite{Meixner}.

\begin{Prop}
\label{Prop:One-variable-Meixner}
Let $\set{P_n(x)}$ be a family of Sheffer polynomials, that is, a polynomial family with a generating function of the form
\begin{equation*}
\sum_{n=0}^\infty \frac{1}{n!} P_n(x) z^n = \exp \Bigl( x U(z) - R(U(z)) \Bigr).
\end{equation*}
Here the conditions on $U$ and $R$ are the same as above. These polynomials are orthogonal if and only if the following two conditions hold:
\begin{enumerate}
\item
$U(z) = (R'(z))^{<-1>}$, and
\item
$R''(z) = 1 + b R'(z) + c (R'(z))^2$.
\end{enumerate}
\end{Prop}

\noindent
In fact, polynomials satisfying the conditions of Proposition~\ref{Prop:One-variable-Meixner} can be listed explicitly. They consist of polynomials orthogonal with respect to the Gaussian, Poisson, gamma, binomial, negative binomial, and continuous binomial (hyperbolic secant) distributions, all important in probability and statistics. It makes sense therefore to look at the polynomials with the generating function~\eqref{One-variable-free-Sheffer}, which we call the free Sheffer polynomials, and in particular at the polynomials satisfying the conditions of Proposition~\ref{Prop:One-variable-free-Meixner}, which we call the free Meixner polynomials. Here the adjective ``free'' refers to their relation to \emph{free probability} \cite{VoiSF}, see \cite{AnsMeixner,AnsAppell} for more details. These polynomials can also be described explicitly, see Theorem 4 of \cite{AnsMeixner}. They include Chebyshev polynomials of the 2nd kind, and other families whose orthogonality measure may include at most two atoms; they belong to the class investigated by Szeg{\"o} and described on pages 26--28 of \cite{Askey-Wilson}. In particular, the semicircular, free Poisson (Marchenko-Pastur) and free binomial distributions are of this type. See Example~\ref{Example:Product-states} for more details.

\br
The parallel between propositions~\ref{Prop:One-variable-free-Meixner} and \ref{Prop:One-variable-Meixner} can be explained by noticing that they are both particular cases of a more general theorem involving the generating function of a specific basic hypergeometric form, see \cite{Al-Salam-Pollaczek} or Theorem 4.8 of \cite{AnsAppell}. Proposition~\ref{Prop:One-variable-free-Meixner} is recovered for $q=0$, while Proposition~\ref{Prop:One-variable-Meixner} is recovered for $q=1$. The general family of orthogonal polynomials involved are the Al-Salam and Chihara polynomials; in particular, the (Rogers) continuous $q$-Hermite polynomials interpolate between the Hermite polynomials and the Chebyshev polynomials of the 2nd kind.

\br
Despite the similarity between single-variable Propositions~\ref{Prop:One-variable-free-Meixner} and \ref{Prop:One-variable-Meixner}, the key point about Definition~\ref{Defn:Free-Sheffer} is that it involves polynomials in non-commuting variables. In contrast, natural multivariate generalizations of Proposition~\ref{Prop:One-variable-Meixner} involve more familiar polynomials in commuting variables, orthogonal with respect to $n$-dimensional measures. They have been investigated by a number of people, see for example \cite{Feinsilver-orthogonal,Letac-Quadratic,Casalis-Wishart,Pommeret-Test,Casalis-Simple-quadratic,Pommeret-Sheffer}. This analysis is usually performed in the context of \emph{natural exponential families}. So this paper may be a precursor to ``free statistics''. For the moment, there are two other motivations for it. First, the hope is that these objects will turn out to play a role in free probability. Second, while there is some work on orthogonal polynomials in non-commuting variables \cite{Constantinescu}, the field appears to be largely unexplored. In particular, while there are many interesting examples of multivariate orthogonal polynomials in commuting variables \cite{Dunkl-Xu}, there is a paucity of examples in the non-commutative case. The original motivation for this paper was to provide such examples. They come from free product states (Section~\ref{Section:Tracial}), and from a certain exponentiation of a free semicircular system (Section~\ref{Section:Example}).

\section{Preliminaries}

\subsection{Polynomials}
Let $\mf{R}\langle \mb{x} \rangle = \mf{R}\langle x_1, x_2, \ldots, x_n \rangle$ be all the polynomials with real coefficients in $n$ non-commuting variables. Multi-indices are elements $\vec{u} \in \set{1, \ldots, n}^k$ for $k \geq 0$; for $\abs{\vec{u}} = 0$ we denote $\vec{u}$ by $\emptyset$. For two multi-indices $\vec{u}, \vec{v}$, denote by $(\vec{u}, \vec{v})$ their concatenation. For $\vec{u}$ with $\abs{\vec{u}} = k$, denote
\[
(\vec{u})^{op} = (u(k), \ldots, u(2), u(1)).
\]
Define an involution on $\mf{R}\langle \mb{x} \rangle$ via an $\mf{R}$-linear extension of
\[
(x_{\vec{u}})^\ast = x_{(\vec{u})^{op}}.
\]
Here $x_{\vec{u}}$ is the monomial $x_{u(1)} \ldots x_{u(k)}$.

\br
A \emph{monic polynomial family} in $\mb{x}$ is a family $\set{P_{\vec{u}}(\mb{x})}$ indexed by all multi-indices
\[
\bigcup_{k=1}^\infty \set{\vec{u} \in \set{1, \ldots, n}^k}
\]
(with $P_{\emptyset} = 1$ being understood) such that $\set{P_{\vec{u}}(\mb{x})} = x_{\vec{u}} + $ \textsl{lower-order terms}. Note that $P_{\vec{u}}^\ast \neq P_{(\vec{u})^{op}}$ in general.

\br
A polynomial family $\set{P_{\vec{u}}}$ is \emph{pseudo-orthogonal} with respect to a functional $\phi$ if
\[
\state{P_{\vec{u}}^\ast P_{\vec{v}}} = 0
\]
whenever $\abs{\vec{u}} \neq \abs{\vec{v}}$ (including $\vec{u} = \emptyset$). The family is \emph{orthogonal} if this is the case whenever $\vec{u} \neq \vec{v}$. Note that an orthogonal polynomial family $\set{P_{\vec{u}}}$ determines its unital functional of orthogonality $\phi$ via $\state{1} = 1$, $\state{P_{\vec{u}}} = 0$, so such a functional, if it exists, is unique.

\br
Most theorems about formal power series remain valid for non-commuting variables. In particular, a series $F(\mb{z}) = 1 + $ \textsl{higher-order terms} has a unique inverse with respect to multiplication, always denoted by $F^{-1}$. Also, an $n$-tuple of series $\mb{U}$ with $U_i(\mb{z}) = z_i + $ \textsl{higher-order terms} has a unique inverse with respect to composition, always denoted by $\mb{U}^{<-1>}$.

\subsection{Linear functionals and free cumulants}
Let $\phi$ be a unital real linear functional on $\mf{R} \langle \mb{x} \rangle$. It can be thought of as a moment functional of the variables $\set{x_1, \ldots, x_n}$. Here their joint moments are
\[
M[x_{\vec{u}}] = \state{x_{\vec{u}}} = \state{x_{u(1)} \ldots x_{u(k)}}.
\]
Denote by
\[
M(\mb{z})
= \sum_{k=1}^\infty \sum_{\abs{\vec{u}} = k} M[x_{\vec{u}}] z_{\vec{u}}
= \state{(1 - \mb{x} \cdot \mb{z})^{-1}} - 1
\]
the ordinary moment generating function of $\phi$. Here, and in the sequel,
\[
\mb{z} = (z_1, z_2, \ldots, z_n)
\]
are formal non-commuting indeterminates, which commute with the $\set{x_i}$. Note that $M(\mb{z})$ completely determines $\phi$.

\br
The free cumulant functional $R$ corresponding to $\phi$ is usually defined using the lattice of non-crossing partitions: $R[1] = 0$ and for $\abs{\vec{u}} = k$,
\begin{equation*}
R[x_{\vec{u}}]
= M[x_{\vec{u}}] - \sum_{\substack{\pi \in \NC(k), \\ \pi \neq \hat{1}}} \prod_{B \in \pi} R \Bigl[\prod_{i \in B} x_{u(i)} \Bigr],
\end{equation*}
which expresses $R[x_{\vec{u}}]$ in terms of the joint moments and sums of products of lower-order free cumulants. From these, we can form the free cumulant generating function via
\begin{equation}
\label{Non-crossing}
\mb{R}(\mb{z})
= \sum_{k=1}^\infty \sum_{\abs{\vec{u}} = k} R[x_{\vec{u}}] z_{\vec{u}}.
\end{equation}
However, in this paper we will not use non-crossing partitions. So for the rest of the paper, we take as the definition of free cumulants the following implicit functional relation, see Section 13 of \cite{Nica-Lecture} or Proposition 3.1 of \cite{AnsAppell}:
\begin{equation}
\label{Cumulant-moment}
\mb{R} \Bigl(w_1 \bigl(1 + \mb{M}(\mb{w}) \bigr), \ldots, w_n \bigl(1 + \mb{M}(\mb{w}) \bigr) \Bigr)
= \mb{M}(\mb{w}).
\end{equation}

\br
To simplify notation, we will assume throughout the paper that the $\set{x_i}$ are centered and have unit covariance,
\[
R[x_i] = \state{x_i} = 0
\]
and
\[
R[x_i x_j] = \state{x_i x_j} - \state{x_i} \state{x_j} = \delta_{ij}.
\]
The results can be modified for more general (in particular, degenerate) covariance, but the formulas become more complicated.

\br
A \emph{state} on $\mf{R} \langle \mb{x} \rangle$ is a linear functional that is unital (that is, $\state{1} = 1$) and positive, that is, for any polynomial $A(\mb{x})$,
\[
\state{A(\mb{x})^\ast A(\mb{x})} \geq 0.
\]
Such a functional cannot always be extended from $\mf{R} \langle \mb{x} \rangle$ to a state on some $C^{\ast}$-algebra. This is already true in the commutative case: a positive functional on $\mf{R}[x_1, x_2]$ need not come from a positive measure on $\mf{R}^2$. The issue is whether the moment problem is solvable; for an example of a non-commutative result, see \cite{Hadwin-Moment}.

\br
A state is \emph{faithful} if $\state{A(\mb{x})^\ast A(\mb{x})} = 0$ only for $A(\mb{x}) = 0$. We will only consider faithful states in this paper; but see Remark~\ref{Remark:Faithful}.

\br
For unital linear functionals
\[
\set{\phi_i \text{ on } \mf{R}[x_i]}_{i=1}^n,
\]
their \emph
{free product functional} $\phi$ on $\mf{R} \langle x_1, x_2, \ldots, x_n \rangle$ is defined by the requirement that
\[
R_\phi[x_{\vec{u}}] = 0
\]
unless all $u(j)$ are equal (that is, mixed free cumulants are zero), and
\[
R_\phi[x_i^k] = R_{\phi_i}[x_i^k].
\]
$\phi$ is a state if $\phi_i$'s are. Conversely, if $\phi$ happens to be of this form, we say that $\set{x_i}$ are freely independent with respect to it. See \cite{VoiSF} or \cite{SpeNCReview} for a lot more about this, and in particular for an explanation of the terminology. If a similar definition is given for the algebra of commutative polynomials in terms of the usual cumulants, one obtains exactly product states, corresponding to product measures, and the notion of independence.

\begin{Ex}
If $\phi_1$ is a state on $\mf{R}[x_1]$, $\phi_2$ is a state on $\mf{R}[x_2]$, and $\phi$ is their free product state on $\mf{R}\langle x_1, x_2 \rangle$, then
\begin{multline*}
\state{x_1 x_2 + x_1 x_2 x_1 + x_1 x_2 x_1 x_2}
= \phi_1[x_1] \phi_2[x_2] + \phi_1[x_1^2] \phi_2[x_2] \\
+ \Bigl( \phi_1[x_1]^2 \phi_2[x_2^2] + \phi_1[x_1^2] \phi[x_2]^2 - \phi_1[x_1]^2 \phi_2[x_2]^2 \Bigr).
\end{multline*}
\end{Ex}

\subsection{Operators}
Define the left partial derivative with respect to $z_i$, $D_i$ by
\begin{gather*}
D_i(1) = 0, \\
D_i z_j = \delta_{ij}, \\
D_i (z_j z_{\vec{u}}) = \delta_{ij} z_{\vec{u}}.
\end{gather*}
Denote by $\mb{D} = (D_1, D_2, \ldots, D_n)$ the left gradient.

\br
Given a monic polynomials family $\set{P_{\vec{u}}}$, define the right partial lowering operator with respect $x_i$, $L_i$, by
\begin{gather*}
L_i (1) = 0, \\
L_i P_j(\mb{x}) = \delta_{ij}, \\
L_i P_{(\vec{u}, j)}(\mb{x}) = \delta_{ij} P_{\vec{u}}(\mb{x}).
\end{gather*}

\section{Free Meixner families}
\begin{Prop}
\label{Prop:Orthogonal-recursion}
Monic polynomials are orthogonal with respect to some faithful state if and only if they satisfy a recursion
\begin{equation}
\label{Recursion}
x_i P_{\vec{u}} = P_{(i, \vec{u})} + \sum_{\abs{\vec{w}} = \abs{\vec{u}}} B_{i, \vec{w}, \vec{u}} P_{\vec{w}} + \sum_{\abs{\vec{v}} = \abs{\vec{u}}-1} C_{i, \vec{v}, \vec{u}} P_{\vec{v}}
\end{equation}
with
\begin{enumerate}
\item
$C_{i, \vec{s}, \vec{u}} = 0$ unless $\vec{u} = (i, \vec{s})$, and $C_{i, \vec{s}, (i, \vec{s})} > 0$,
\item
Denoting $\vec{s}_j = (s(j), \ldots, s(k))$,
\[
B_{i, \vec{s}, \vec{u}} \prod_{j=1}^k C_{s(j), \vec{s}_{j+1}, \vec{s}_j} = B_{i, \vec{u}, \vec{s}} \prod_{j=1}^k C_{u(j), \vec{u}_{j+1}, \vec{u}_j}.
\]
\end{enumerate}
\end{Prop}

\begin{proof}
First assume that the polynomials are orthogonal with respect to some faithful state $\phi$. Denote
\[
\ip{S(\mb{x})}{T(\mb{x})} = \state{S(\mb{x})^\ast T(\mb{x})}.
\]
Since the polynomials are monic, for any fixed $\vec{u}$, $i$,
\[
x_i P_{\vec{u}} = P_{(i, \vec{u})} + \sum_{\abs{\vec{v}} \leq \abs{\vec{u}}} \alpha_{i, \vec{v}, \vec{u}} P_{\vec{v}}
\]
for some coefficients $ \alpha_{i, \vec{v}, \vec{u}}$. Also,
\[
\ip{P_{\vec{v}}}{x_i P_{\vec{u}}}
= \state{P_{\vec{v}}^\ast (x_i P_{\vec{u}})}
= \state{(x_i P_{\vec{v}})^\ast P_{\vec{u}}}
= \ip{x_i P_{\vec{v}}}{P_{\vec{u}}}
= 0
\]
for $\abs{\vec{v}} \leq \abs{\vec{u}} - 2$. It follows that the polynomials satisfy a recursion of the type \eqref{Recursion}. In that case, for general $\vec{u}$ and $\vec{s}$
\begin{equation}
\label{Associativity}
\begin{split}
\ip{x_i P_{\vec{u}}}{P_{\vec{s}}}
& = \ip{P_{(i, \vec{u})}}{P_{\vec{s}}}
+ \sum_{\abs{\vec{w}} = \abs{\vec{u}}} B_{i, \vec{w}, \vec{u}} \ip{P_{\vec{w}}}{P_{\vec{s}}}
+ \sum_{\abs{\vec{v}} = \abs{\vec{u}} - 1} C_{i, \vec{v}, \vec{u}} \ip{P_{\vec{v}}}{P_{\vec{s}}} \\
 = \ip{P_{\vec{u}}}{x_i P_{\vec{s}}}
& = \ip{P_{\vec{u}}}{P_{(i, \vec{s})}}
+ \sum_{\abs{\vec{w}} = \abs{\vec{s}}} B_{i, \vec{w}, \vec{s}} \ip{P_{\vec{u}}}{P_{\vec{w}}}
+ \sum_{\abs{\vec{v}} = \abs{\vec{s}} - 1} C_{i, \vec{v}, \vec{s}} \ip{P_{\vec{u}}}{P_{\vec{v}}}
\end{split}
\end{equation}
Pseudo-orthogonality implies that for $\abs{\vec{s}} = \abs{\vec{u}} -1$
\[
\sum_{\abs{\vec{v}} = \abs{\vec{u}} - 1} C_{i, \vec{v}, \vec{u}} \ip{P_{\vec{v}}}{P_{\vec{s}}} = \ip{P_{\vec{u}}}{P_{(i, \vec{s})}},
\]
and for $\abs{\vec{s}} = \abs{\vec{u}}$
\[
\sum_{\abs{\vec{w}} = \abs{\vec{u}}} B_{i, \vec{w}, \vec{u}} \ip{P_{\vec{w}}}{P_{\vec{s}}} = \sum_{\abs{\vec{w}} = \abs{\vec{s}}} B_{i, \vec{w}, \vec{s}} \ip{P_{\vec{u}}}{P_{\vec{w}}}.
\]
(the case $\abs{\vec{s}} = \abs{\vec{u}} + 1$ is redundant). Using the orthogonality assumption,
\[
C_{i, \vec{s}, \vec{u}} V_{\vec{s} \vec{s}} = \delta_{\vec{u}, (i, \vec{s})} V_{\vec{u} \vec{u}}
\]
and
\[
B_{i, \vec{s}, \vec{u}} V_{\vec{s} \vec{s}} = B_{i, \vec{u}, \vec{s}} V_{\vec{u} \vec{u}},
\]
where
\[
V_{\vec{u} \vec{u}}
= \ip{P_{\vec{u}}}{P_{\vec{u}}}.
\]
It follows that
\begin{equation}
\label{V_uu}
V_{\vec{u} \vec{u}} = \prod_{j=1}^k C_{\vec{u}_{j+1}, \vec{u}_j, u(j)},
\end{equation}
\begin{equation}
\label{C}
C_{i, \vec{s}, \vec{u}} \prod_{j=1}^{k-1} C_{s(j), \vec{s}_{j+1}, \vec{s}_j} = \delta_{\vec{u}, (i, \vec{s})} \prod_{j=1}^k C_{u(j), \vec{u}_{j+1}, \vec{u}_j}
\end{equation}
and
\begin{equation}
\label{B}
B_{i, \vec{s}, \vec{u}} \prod_{j=1}^k C_{s(j), \vec{s}_{j+1}, \vec{s}_j} = B_{i, \vec{u}, \vec{s}} \prod_{j=1}^k C_{u(j), \vec{u}_{j+1}, \vec{u}_j}.
\end{equation}
Equation~\eqref{B} is condition (b). Equation~\eqref{C} is equivalent to requiring that $C_{i, \vec{s}, \vec{u}} = 0$ unless $\vec{u} = (i, \vec{s})$, and faithfulness of $\phi$ implies that $C_{i, \vec{s}, (i, \vec{s})} > 0$, which together form condition (a).

\br
Conversely, assume that the polynomials satisfy the recursion~\eqref{Recursion} with the conditions of the proposition. On $\mf{R} \langle \mb{x} \rangle$, define the functional $\varphi$ by requiring that the induced inner product
\[
\ip{S(\mb{x})}{T(\mb{x})} = \state{S(\mb{x})^\ast T(\mb{x})}
\]
satisfies
\[
\ip{P_{\vec{u}}}{P_{\vec{v}}}
= \delta_{\vec{u} \vec{v}} V_{\vec{u} \vec{u}}
= \begin{cases}
0 & \text{if } \vec{u} \neq \vec{v}, \\
V_{\vec{u} \vec{u}} & \text{if } \vec{u} = \vec{v},
\end{cases}
\]
where $V_{\vec{u} \vec{u}}$ is now \emph{defined} via equation~\eqref{V_uu}, and extending linearly. So for $S(\mb{x}) = \sigma_{\emptyset} + \sum_{\vec{u}} \sigma_{\vec{u}} P_{\vec{u}}(\mb{x})$, $T(\mb{x}) = \tau_{\emptyset} + \sum_{\vec{u}} \tau_{\vec{u}} P_{\vec{u}}(\mb{x})$,
\[
\ip{S(\mb{x})}{T(\mb{x})} = \sigma_{\emptyset} \tau_{\emptyset} + \sum_{\vec{u}} \sigma_{\vec{u}} \tau_{\vec{u}} V_{\vec{u} \vec{u}}.
\]
If this functional is well-defined, the given polynomials are orthogonal with respect to it. Also, since $V_{\vec{u} \vec{u}}$ are positive, the functional will be positive and faithful.

\br
To show that this definition is consistent, we need to show that if
\begin{equation}
\label{ST}
S(\mb{x}) T(\mb{x}) = S'(\mb{x}) T'(\mb{x}),
\end{equation}
then
\[
\ip{S^\ast}{T} = \ip{(S')^\ast}{T'}.
\]
For $\mf{R} \langle \mb{x} \rangle$, the fundamental theorem of algebra no longer holds, but these polynomials still form a Unique Factorization Domain. Thus the equality \eqref{ST} reduces to the situation $(Q S) T = Q (S T)$. By linearity, we may assume that $S$ is a monomial. But in that case, by iteration we may assume that $S = x_i$. Finally, by linearity again we may assume that $Q^\ast, T$ are basis polynomials. Thus we only need to satisfy the following condition:
\[
\begin{split}
\ip{x_i P_{\vec{u}}}{P_{\vec{s}}}
& = \ip{P_{\vec{u}}}{P_{(i, \vec{s})}}
\end{split}
\]
which, using the recursion relation, is equivalent to equation~\eqref{Associativity}. The arguments from the first half of the proof imply that this equality holds provided that conditions (a), (b) are satisfied.
\end{proof}

\begin{Remark}
It follows from the proof of the preceding proposition that any pseudo-or\-tho\-go\-nal polynomials satisfy a recursion of type \eqref{Recursion}.
\end{Remark}

\begin{Remark}
\label{Remark:Faithful}
If the appropriate part of condition (a) of the proposition is replaced by the condition $C_{i, \vec{s}, (i, \vec{s})} \geq 0$, it follows that the corresponding polynomials are still orthogonal with respect to a state that need not be faithful. The converse characterization is an interesting question that is not treated in this paper.
\end{Remark}

\begin{Lemma}
\label{Lemma:Resolvent}
Let $\set{P_{\vec{u}}}$ be a family of free Sheffer polynomials as in Definition~\ref{Defn:Free-Sheffer}, with
\[
H(\mb{x}, \mb{z}) = 1 + \sum_{\vec{u}} P_{\vec{u}}(\mb{x}) z_{\vec{u}} = F(\mb{z}) \Bigl( 1 - \mb{x} \cdot \mb{V}(\mb{z}) \Bigr)^{-1}.
\]
Assume more particularly that $F(\mb{z}) = 1 - \sum_{i=1}^n z_i^2 + $ \textsl{higher-order terms}. Define the functional $\phi$ on $\mf{R}\langle x_1, x_2, \ldots, x_n \rangle$ by $\state{1} = 1$, $\state{P_{\vec{u}}} = 0$ for $\abs{\vec{u}} \geq 1$. Then in fact,
\[
H(\mb{x}, \mb{z}) = \Bigl( 1 - \mb{x} \cdot \mb{U}(\mb{z}) + R(\mb{U}(\mb{z})) \Bigr)^{-1},
\]
where $R(\mb{z})$ is the free cumulant generating function of $\phi$, and $U_i(\mb{z}) = V_i(\mb{z}) F^{-1}(\mb{z})$. We say that $\set{P_{\vec{u}}}$ is the free Sheffer family associated to the functional $\phi$ and the functions $\mb{U}$. Note that if a free Sheffer family is orthogonal, it is orthogonal with respect to the functional $\varphi$ to which it is associated.
\end{Lemma}

\begin{proof}
By definition of $\phi$ and $H$, $\state{H(\mb{x}, \mb{z})} = 1$. Then
\[
1 = F(\mb{z}) \state{\bigl( 1 - \mb{x} \cdot \mb{U}(\mb{z}) F(\mb{z}) \bigr)^{-1}} = F(\mb{z}) (1 + M(\mb{U}(\mb{z}) F(\mb{z})))
\]
Since the $n$-tuple of power series $\mb{U}$ is invertible under composition, we may write
\[
F(\mb{z}) = \bigl(1 + K(\mb{U}(\mb{z})) \bigr)^{-1}
\]
for some power series $K$. Then
\[
1 + K(\mb{U}) = 1 + M \Bigl( \mb{U} \bigl(1 + K(\mb{U}) \bigr)^{-1} \Bigr).
\]
Therefore from equation~\eqref{Cumulant-moment},
\[
\begin{split}
& R \Bigl( U_1 \bigl(1 + K(\mb{U}) \bigr)^{-1} \Bigl( 1 + M \bigl( \mb{U} \bigl(1 + K(\mb{U}) \bigr)^{-1} \bigr) \Bigr), \ldots, \\
&\qquad \ldots, U_n \bigl(1 + K(\mb{U}) \bigr)^{-1} \Bigl( 1 + M \bigl( \mb{U} \bigl(1 + K(\mb{U}) \bigr)^{-1} \bigr) \Bigr) \Bigr) = K(\mb{U}).
\end{split}
\]
However, this expression also equals
\[
R \Bigl( U_1  \bigl(1 + K(\mb{U}) \bigr)^{-1} (1 + K(\mb{U})), \ldots, U_n  \bigl(1 + K(\mb{U}) \bigr)^{-1} (1 + K(\mb{U})) \Bigr)
= R(\mb{U}).
\]
Thus $F(\mb{z}) = \bigl(1 + R(\mb{U}(\mb{z})) \bigr)^{-1}$ and
\[
F(\mb{z}) \Bigl( 1 - \mb{x} \cdot \mb{V}(\mb{z}) \Bigr)^{-1} = \Bigl( 1 - \mb{x} \cdot \mb{U}(\mb{z}) + R(\mb{U}(\mb{z})) \Bigr)^{-1}.
\]
\end{proof}

\begin{Prop}\cite[Theorem 3.21]{AnsAppell}
\label{Prop:Inverse}
Suppose that a family of free Sheffer polynomials is pseudo-orthogonal. Then for $R$, $\mb{U}$ as in Lemma~\ref{Lemma:Resolvent},
\[
(D_i R) (\mb{U}(\mb{z})) = z_i.
\]
\end{Prop}

\begin{Remark}
Both $\mb{D} R$ and $\mb{U}$ are $n$-tuples of non-commutative power series invertible under composition. So
\begin{enumerate}
\item
Given $R$, the preceding proposition completely determines $\mb{U}$, and vice versa. From now on, we will always assume this relationship between $R$ and $\mb{U}$.
\item
Since the inverse under composition is unique, also
\begin{equation}
\label{Inversion2}
U_i((\mb{D} R)(\mb{z})) = z_i.
\end{equation}
\end{enumerate}
\end{Remark}

\begin{Defn}
A state $\phi$ on $\mf{R} \langle \mb{x} \rangle$ is called a \emph{free Meixner state} if, for $R$ its free cumulant generating function and $\mb{U}$ determined by the preceding remark, the free Sheffer polynomials with the generating function
\[
\Bigl( 1 - \mb{x} \cdot \mb{U}(\mb{z}) + R(\mb{U}(\mb{z})) \Bigr)^{-1}
\]
are orthogonal.
\end{Defn}

\begin{Thm}
\label{Thm:Quadratic}
Suppose that a family of free Sheffer polynomials is pseudo-orthogonal. Then
\begin{enumerate}
\item
The power series $\mb{U}$ satisfy the relation
\[
z_j = U_j + \sum_{i,t} B_{ij}^{t} U_i z_t + \sum_{i,s,t} C_{ij}^{st} U_i z_s z_t .
\]
In other words, denoting by $A$ the matrix
\[
I + \sum_t B^t z_t + \sum_{s,t} C^{st} z_s z_t,
\]
$\mb{U} =  \mb{z} A^{-1}$.
\item
The polynomials satisfy the recursion
\begin{gather*}
x_i P_{s} = P_{(i,s)} + \sum_j B_{ij}^{s} P_{j} + \delta_{is}, \\
x_i P_{(s,t,\vec{u})} = P_{(i,s,t,\vec{u})} + \sum_j B_{ij}^{s} P_{(j,t,\vec{u})} + \sum_j (\delta_{is} \delta_{jt} + C_{ij}^{st}) P_{(j, \vec{u})}.
\end{gather*}
\item
The free cumulant generating function satisfies
\[
D_i D_j R(\mb{z}) = \delta_{ij} + \sum_t B_{ij}^{t} D_t R(\mb{z}) + \sum_{s,t} C_{ij}^{st} D_s R(\mb{z}) \ D_t R(\mb{z}).
\]
\end{enumerate}
\end{Thm}

\begin{proof}
By definition of the function $H$ in Lemma~\ref{Lemma:Resolvent},
\[
L_j H(\mb{x}, \mb{z}) = H(\mb{x}, \mb{z}) z_j.
\]
Also from that lemma,
\[
(1 + R(\mb{U}(\mb{z}))) H = (\mb{x} \cdot \mb{U}(\mb{z})) H + 1.
\]
Applying $L_j$ to this expression, we get
\[
L_j \bigl( (x \cdot \mb{U}(\mb{z})) H \bigr)
= L_j \bigl( (1 + R(\mb{U}(\mb{z}))) H \bigr)
= (1 + R(\mb{U}(\mb{z}))) H z_j
= (\mb{x} \cdot \mb{U}(\mb{z})) H z_j + z_j.
\]
Expanding $H$ in powers of $z$, we get
\[
L_j \Bigl( \sum_i x_i U_i(z) (1 + \sum_{\vec{u}} P_{\vec{u}} z_{\vec{u}}) \Bigr) = \sum_i x_i U_i(z) (1 + \sum_{\vec{u}} P_{\vec{u}} z_{\vec{u}}) z_j + z_j,
\]
and so
\begin{equation}
\label{Expansion1}
U_j + \sum_{\vec{u}, i} L_j(x_i P_{\vec{u}}) U_i(z) z_{\vec{u}} = z_j + \sum_i x_i U_i(z) z_j + \sum_{\vec{u}, i} x_i P_{\vec{u}} U_i(z) z_{\vec{u}} z_j,
\end{equation}
where we used the fact that $L_j(x_i) = L_j(P_i) = \delta_{ij}$.

\br
Since $U_i = z_i + $ \textsl{higher-order terms},
\begin{equation}
\label{z-U-expansion}
z_j = U_j + \sum_{i, \vec{u}} a_{i, j, \vec{u}} U_i(z) z_{\vec{u}}
\end{equation}
for some coefficients $\set{a_{i, j, \vec{u}}}$. Using equation~\eqref{Inversion2},
\[
D_j R = z_j + \sum_{i, \vec{u}} a_{i, j, \vec{u}} z_i (\mb{D} R)_{\vec{u}},
\]
where $(\mb{D} R)_{\vec{u}} = (D_{u(1)} R) (D_{u(2)} R) \ldots (D_{u(k)} R)$. Therefore
\begin{equation}
\label{DDR}
D_i D_j R = \delta_{ij} + \sum_{\vec{u}} a_{i, j, \vec{u}} (\mb{D} R)_{\vec{u}}.
\end{equation}

\br
Combining equations~\eqref{Expansion1} and \eqref{z-U-expansion},
\[
\begin{split}
U_j + \sum_{i, \vec{u}} L_j(x_i P_{\vec{u}}) U_i(z) z_{\vec{u}} = U_j
& + \sum_{i, \vec{u}} a_{i, j, \vec{u}} U_i(z) z_{\vec{u}} \\
& + \sum_i x_i U_i(z) z_j + \sum_{i, \vec{u}} x_i P_{\vec{u}} U_i(z) z_{\vec{u}} z_j.
\end{split}
\]
Equating coefficients of $U_i z_{\vec{u}}$,
\[
L_j(x_i P_{\vec{u}}) = a_{i, j, \vec{u}} + \delta_{u(k), j} x_i P_{(u(1), \ldots, u(k-1))}.
\]
Since the polynomials are pseudo-orthogonal, they satisfy a recursion relation~\eqref{Recursion}. So
\[
L_j (x_i P_{\vec{u}}) = \delta_{u(k), j} P_{(i, u(1), \ldots, u(k-1))} + \sum_{\abs{\vec{w}} = \abs{\vec{u}}-1} B_{i, (\vec{w}, j), \vec{u}} P_{\vec{w}} + \sum_{\abs{\vec{v}} = \abs{\vec{u}}-2} C_{i, (\vec{v}, j), \vec{u}} P_{\vec{v}}.
\]
Combining the two preceding equations with equation~\eqref{Recursion} for $x_i P_{(u(1), \ldots, u(k-1))}$, we get
\begin{gather*}
\begin{split}
a_{i, j, \vec{u}} + \delta_{u(k), j} \Bigl( P_{(i, u(1), \ldots, u(k-1))}
& + \sum_{\abs{\vec{w}} = \abs{\vec{u}}-1} B_{i, \vec{w}, (u(1), \ldots, u(k-1))} P_{\vec{w}} \\
& + \sum_{\abs{\vec{v}} = \abs{\vec{u}}-2} C_{i, \vec{v}, (u(1), \ldots, u(k-1))} P_{\vec{v}} \Bigr)
\end{split} \\
= \delta_{u(k), j} P_{(i, u(1), \ldots, u(k-1))} + \sum_{\abs{\vec{w}} = \abs{\vec{u}}-1} B_{i, (\vec{w}, j), \vec{u}} P_{\vec{w}} + \sum_{\abs{\vec{v}} = \abs{\vec{u}}-2} C_{i, (\vec{v}, j), \vec{u}} P_{\vec{v}}.
\end{gather*}
Equating coefficients,
\[
\begin{split}
a_{i, j, \vec{u}}
& = \sum_{\abs{\vec{w}} = \abs{\vec{u}}-1} \left( B_{i, (\vec{w}, j), \vec{u}} - \delta_{u(k), j} B_{i, \vec{w}, (u(1), \ldots, u(k-1))} \right) P_{\vec{w}} \\
&\quad + \sum_{\abs{\vec{v}} = \abs{\vec{u}}-2} \left( C_{i, (\vec{v}, j), \vec{u}} - \delta_{u(k), j} C_{i, \vec{v}, (u(1), \ldots, u(k-1))} \right) P_{\vec{v}}.
\end{split}
\]
In particular, for $\vec{u} = t$ this says
\[
a_{i,j,t} = B_{i,j,t} - \delta_{jt} B_{i, \emptyset, \emptyset},
\]
and for $\vec{u} = (s,t)$ this says
\[
a_{i, j, (s,t)}
= \sum_{w} \left( B_{i, (w,j), (s,t)} - \delta_{jt} B_{i, w, s} \right) P_{w} + \left( C_{i, j, (s,t)} - \delta_{jt} C_{i, \emptyset, s} \right).
\]
Therefore
\[
B_{i, (\vec{w}, j), \vec{u}} = \delta_{u(k), j} B_{i, \vec{w}, (u(1), \ldots, u(k-1))}; \qquad B_{i, j, t} = \delta_{jt} B_{i, \emptyset, \emptyset} + a_{i,j,t},
\]
\[
C_{i, (\vec{v}, j), \vec{u}} = \delta_{u(k), j} C_{i, \vec{v}, (u(1), \ldots, u(k-1))}; \qquad C_{i, j, (s,t)} = \delta_{jt} C_{i, \emptyset, s} + a_{i, j, (s,t)}.
\]
So
\[
B_{i, (j, \vec{u}), (t, \vec{u})} = \delta_{jt} B_{i, \emptyset, \emptyset} + a_{i,j,t},
\]
\[
C_{i, (j, \vec{u}), (s, t, \vec{u})} = \delta_{jt} C_{i, \emptyset, s} + a_{i, j, (s,t)}
\]
and zero otherwise.
\[
x_i = P_i + B_{i, \emptyset, \emptyset},
\]
\[
x_i P_t = P_{(i,t)} + \sum_s B_{i, s, t} P_s + C_{i, \emptyset, t}.
\]
So $B_{i, \emptyset, \emptyset} = R[x_i] = 0$, $C_{i, \emptyset, t} = R[x_i  x_t] = \delta_{it}$. Also, $a_{i, j, \vec{u}} = 0$ for $\abs{\vec{u}} > 2$. Denote $B_{ij}^{t} = a_{i,j,t}$, $C_{ij}^{st} = a_{i, j, (s,t)}$. Part (b) follows. For parts (a) and (c), use equations \eqref{z-U-expansion}, \eqref{DDR}, respectively.
\end{proof}

\begin{Cor}
\label{Cor:Orthogonal}
Let $\phi$ be a state, $R$ its free cumulant generating function, $\mb{U}$ the corresponding power series determined by Proposition~\ref{Prop:Inverse}, and $\set{P_{\vec{u}}}$ the corresponding free Sheffer polynomials. $\phi$ is a faithful free Meixner state if and only if the following equivalent conditions hold:
\[
z_j = U_j + \sum_{i,t} B_{ij}^{t} U_i z_t + \sum_i C_{ij} U_i z_i z_j,
\]
or
\begin{equation}
\label{Quadratic-orthogonal}
D_i D_j R = \delta_{ij} + \sum_t B_{ij}^{t} D_t R + C_{ij} D_i R \ D_j R,
\end{equation}
or
\begin{gather}
x_i P_{t} = P_{(i,t)} + \sum_j B_{ij}^{t} P_{j} + \delta_{it}, \notag \\
x_i P_{(t, \vec{u})} = P_{(i, t, \vec{u})} + \sum_j B_{ij}^{t} P_{(j, \vec{u})} + \delta_{it} (1 + C_{i, u(1)}) P_{\vec{u}}. \label{corollary-recursion}
\end{gather}
In all cases, the coefficients have to satisfy
\begin{enumerate}
\item
$C_{ij} > -1$.
\item
$B_{ij}^{t} = B_{it}^{j}$.
\item
For each $j,t$, either $B_{ij}^{t} = 0$ for all $i$, or $C_{ju} = C_{tu}$ for all $u$.
\end{enumerate}
\end{Cor}

\begin{proof}
If the free Sheffer polynomials are orthogonal with respect to the state $\varphi$, then in particular $\state{P_{\vec{u}}} = 0$, so by Lemma~\ref{Lemma:Resolvent}, $\varphi$ is exactly the state with the free cumulant generating function $R$.

\br
Combine Proposition~\ref{Prop:Orthogonal-recursion} with Theorem~\ref{Thm:Quadratic}. It follows that $C_{st}^{ij} = \delta_{is} \delta_{jt} C_{ij}$ and
\[
C_{i, \vec{u}, (i, \vec{u})} = 1 + C_{i, u(1)}
\]
and zero otherwise, so condition (a) follows from Proposition~\ref{Prop:Orthogonal-recursion}(a). Also,
\[
B_{i, (j, \vec{u}), (t, \vec{u})} V_{(j,\vec{u}), (j,\vec{u})} = B_{i, (t, \vec{u}), (j, \vec{u})} V_{(t,\vec{u}), (t,\vec{u})},
\]
so
\[
B_{ij}^{t} V_{(j,\vec{u}), (j,\vec{u})} = B_{it}^{j} V_{(t,\vec{u}), (t,\vec{u})}.
\]
For $\abs{\vec{u}} = 0$, this says
\[
B_{ij}^{t} = B_{it}^{j},
\]
implying condition (b). For longer $\vec{u}$, this says
\[
B_{ij}^{t} C_{j, \vec{u}, (j, \vec{u})} = B_{it}^{j} C_{t, \vec{u}, (t, \vec{u})},
\]
so
\[
B_{ij}^{t} (1 + C_{ju}) = B_{it}^{j} (1 + C_{tu}),
\]
implying condition (c).

\br
Conversely, suppose that for the state $\varphi$ and the corresponding free Sheffer polynomials $\set{P_{\vec{u}}}$, the recursion~\eqref{corollary-recursion} with conditions (a-c) holds. Then by Proposition~\ref{Prop:Orthogonal-recursion}, the polynomials are orthogonal, necessarily with respect to $\phi$, and $\phi$ is faithful. The equivalence of the conditions for $R$, $\mb{U}$, and the polynomials in the corollary follows from Theorem~\ref{Thm:Quadratic}.

\end{proof}

\section{First-order, tracial case}
\label{Section:Tracial}
\noindent
Throughout this section, we will assume that the state $\phi$ is tracial, that is, for any $S, T$,
\[
\state{S(\mb{x}) T(\mb{x})} = \state{T(\mb{x}) S(\mb{x})}.
\]
This produces two simplifications. First, for any $\vec{u}, i$,
\begin{equation}
\label{Cyclic-R}
R[x_{\vec{u}} x_i] =  R[x_i x_{\vec{u}}].
\end{equation}
This is not apparent from the definition of $R$ via equation~\eqref{Cumulant-moment}, but follows easily from the definition using non-crossing partitions.

\br
Second, any pseudo-orthogonal polynomials can be orthogonalized (with real coefficients).

\begin{Remark}
Starting with an arbitrary monic polynomial family, by using the Gram-Schmidt procedure it can be transformed into a pseudo-orthogonal family; note that this family is still monic. Given an ordering of the monomials of the same degree, the procedure can be applied further to produce an orthogonal family. However, this will necessarily destroy the monic condition. Therefore, the condition that monic orthogonal polynomials exist is rather strong, and does not hold for all tracial states.
\end{Remark}

\begin{Lemma}
\label{Lemma:Cyclic}
Let $B, C$ be as in Theorem~\ref{Thm:Quadratic}. $B_{ij}^t$ is invariant under cyclic permutations of $(j,i,t)$, and $\sum_t B_{ij}^t B_{ct}^d + C_{ij}^{cd}$ is invariant under cyclic permutations of $(j,i,c,d)$.
\end{Lemma}

\begin{proof}
By assumption,
\[
D_i D_j R = \delta_{ij} + \sum_t B_{ij}^t D_t R + \sum_{s,t} C_{ij}^{st} D_s R \ D_t R,
\]
and also
\[
R = \sum_i z_i^2 + \sum R[x_a x_b x_c] z_a z_b z_c + \sum R[x_a x_b x_c x_d] z_a z_b z_c z_d + \ldots.
\]
Then
\[
D_j R = z_j + \sum R[x_j x_b x_c] z_b z_c + \sum R[x_j x_b x_c x_d] z_b z_c z_d + \ldots
\]
and so
\[
\begin{split}
D_i D_j R
& = \delta_{ij} + \sum R[x_j x_i x_c] z_c + \sum R[x_j x_i x_c x_d] z_c z_d + \ldots \\
& = \delta_{ij} + \sum B_{ij}^t z_t + \sum B_{ij}^t R[x_t x_b x_c] z_b z_c + \sum C_{ij}^{cd} z_c z_d + \ldots.
\end{split}
\]
It follows that
\[
R[x_j x_i x_t] = B_{ij}^t
\]
and
\[
R[x_j x_i x_c x_d] = \sum_t B_{ij}^t R[x_t x_c x_d] + C_{ij}^{cd} = \sum_t B_{ij}^t B_{ct}^d + C_{ij}^{cd}.
\]
So the result follows from cyclic symmetry~\eqref{Cyclic-R}.
\end{proof}

\begin{Lemma}
\label{Lemma:Change}
Let $O$ be an orthogonal transformation on $\mf{R}^n$. Perform changes of variables $\mb{x} = O \mb{y}$, $\mb{w} = O^{-1} \mb{z}$. Then
\begin{gather*}
R_{\mb{x}}(\mb{z}) = R_{\mb{y}}(\mb{w}), \\
\mb{D} R_{\mb{x}}(\mb{z}) = O \mb{D} R_{\mb{y}}(\mb{w}), \\
Q R_{\mb{x}}(\mb{z}) = O (Q R_{\mb{y}}(\mb{w})) O^{-1},
\end{gather*}
where $Q_{ij} R = D_i D_j R$.
\begin{equation}
\label{Change-of-variable}
\Bigl( 1 - \mb{x} \cdot \mb{U}(\mb{z}) + R_{\mb{x}}(\mb{U}(\mb{z})) \Bigr)^{-1}
= \Bigl( 1 - \mb{y} \cdot \mb{V}(\mb{w}) + R_{\mb{y}}(\mb{V}(\mb{w})) \Bigr)^{-1}
\end{equation}
for $\mb{V}(\mb{w}) = O^{-1} \mb{U}(O \mb{w})$. The induced functional on $\mf{R} \langle \mb{y} \rangle$ is tracial if $\phi$ was.  The polynomials with the generating function \eqref{Change-of-variable} are orthogonal for all such changes of variable $O$ if and only if, in addition to the conditions of Corollary~\ref{Cor:Orthogonal}, $C_{ij} \equiv c$.
\end{Lemma}

\begin{proof}
\[
R_{\mb{x}}(\mb{z}) = \sum_{\vec{u}} R[x_{\vec{u}}] z_{\vec{u}}.
\]
So by linearity of $R$, for $x_i = \sum_j O_{ij} y_j$,
\[
\begin{split}
R_{\mb{x}}(\mb{z})
= \sum_{\vec{u}} \sum_{\abs{\vec{v}} = \abs{\vec{u}}} \prod_{i=1}^k O_{u(i) v(i)} R[y_{\vec{v}}] z_{\vec{u}}
& = \sum_{\vec{v}} R[y_{\vec{v}}] \sum_{\abs{\vec{v}} = \abs{\vec{u}}} \prod_{i=1}^k O_{u(i) v(i)} z_{\vec{u}} \\
& = \sum_{\vec{v}} R[y_{\vec{v}}] w_{\vec{v}}
= R_{\mb{y}}(\mb{w}).
\end{split}
\]
where $w_j = \sum_i O_{ij} z_i$. Also,
\[
\begin{split}
D_i R_{\mb{x}}(\mb{z})
& = \sum_{\vec{v}} R[y_{\vec{v}}] \sum_{\abs{\vec{v}} = \abs{\vec{u}} + 1} O_{i v(1)} \prod_{j=2}^{k} O_{u(j) v(j)} z_{\vec{u}} \\
& = \sum_{\vec{v}, s} R[y_s y_{\vec{v}}] O_{i s} w_{\vec{v}}
= \sum_s O_{i s} D_s R_{\mb{y}}(\mb{w}).
\end{split}
\]
Similarly,
\[
D_i D_j R_{\mb{x}}(\mb{z})
= \sum_{s,t} O_{is} O_{jt} D_s D_t R_{\mb{y}}(\mb{w}).
\]
Equation~\eqref{Change-of-variable} follows. The tracial property is clear.

\br
If
\[
D_i D_j R_{\mb{x}}(\mb{z}) = \delta_{ij} + \sum_t B_{ij}^t D_t R_{\mb{x}}(\mb{z}) + C_{ij} D_i R_{\mb{x}}(\mb{z}) \ D_j R_{\mb{x}}(\mb{z}),
\]
then
\[
\begin{split}
D_i D_j R_{\mb{y}}(\mb{w}) = \delta_{ij} & + \sum B_{\alpha \beta}^t O_{\alpha i} O_{\beta j} O_{ts} D_s R_{\mb{y}}(\mb{w}) \\
& + \sum C_{s t} O_{\alpha i} O_{\beta j} O_{s u} O_{t v} D_u R_{\mb{y}}(\mb{w}) \ D_v R_{\mb{y}}(\mb{w}).
\end{split}
\]
For orthogonality of the induced free Sheffer polynomials in $\mb{y}$, we check the conditions of Corollary~\ref{Cor:Orthogonal}.  By Lemma~\ref{Lemma:Cyclic},
\[
B_{\alpha \beta}^t O_{\alpha i} O_{\beta j} O_{ts}
= B_{\alpha t}^\beta O_{\alpha i} O_{\beta j} O_{ts}
= B_{\alpha \beta}^t O_{\alpha i} O_{tj} O_{\beta s},
\]
so this expression is symmetric in $j, s$. On the other hand, we also need
\[
\sum_{s, t} C_{st} O_{s i} O_{t j} O_{s u} O_{t v} = \delta_{iu} \delta_{jv} E_{ij}.
\]
Taking the sum of these expressions with respect to $\sum_{i, j} O_{a i} O_{b j}$, we get
\[
C_{ab} O_{a u} O_{b v} =  E_{uv} O_{a u} O_{b v}.
\]
It follows that for all $a,b,u,v$, $C_{ab} = E_{uv}$, hence $C_{ab} \equiv c$.  Finally, for constant $C$ the last condition of the corollary is trivially true.
\end{proof}

\begin{Cor}
\label{Cor:Commute}
Let $B, C$ be as in Corollary~\ref{Cor:Orthogonal}. Then $B_{ij}^t$ is symmetric under all permutations of $(i,j,t)$, and $C_{ij}$ is symmetric in its arguments. If $C \equiv 0$, all the matrices $\set{B^t}$ commute.
\end{Cor}

\begin{proof}
The symmetry of $B_{ij}^t$ comes by combining the cyclic symmetry from Lemma~\ref{Lemma:Cyclic} with the transposition symmetry from Corollary~\ref{Cor:Orthogonal}. Also from that lemma,
\[
R[x_j x_i x_i x_j] = \sum_t B_{ij}^t B_{it}^j + C_{ij},
\]
while
\begin{equation*}
R[x_i x_j x_j x_i] = \sum_t B_{ji}^t B_{jt}^i + C_{ji}.
\end{equation*}
It follows that $C_{ji} = C_{ij}$. Using the cyclic symmetry from the lemma again and setting $C \equiv 0$,
\[
\sum_t B_{ij}^t B_{ct}^d = \sum_t B_{ci}^t B_{dt}^j.
\]
So
\[
(B^j B^c)_{id} = \sum_t B_{it}^j B_{td}^c = \sum_t B_{it}^c B_{td}^j = (B^c B^j)_{id}.
\]
\end{proof}

\begin{Ex}[Product states]
\label{Example:Product-states}
Let $\phi^{b,c}$ be a one-dimensional free Meixner state, that is, the state on $\mf{R}[x]$ whose free cumulant generating function satisfies the equation in Proposition~\ref{Prop:One-variable-free-Meixner}. The solution of this equation is
\[
R(z) = \frac{z^{-1} - b - \sqrt{(z^{-1} - b)^2 - 4c}}{2c}.
\]
Note that the free cumulant generating function $R$ differs from a more familiar $R$-transform by a factor of $z$. $\phi^{b,c}$ itself can be identified with the measure
\[
\frac{1}{2 \pi} \frac{\sqrt{4(1+c) - (x - b)^2}}{1 + b x + c x^2}\,dx + \text{ zero, one, or two atoms};
\]
see Theorem 4 of \cite{AnsMeixner} for a more detailed description, with different normalizations. Here $b \in \mf{R}$, and $c \geq -1$ (for $c = -1$, the measure is purely atomic, so the corresponding state is not faithful). In particular, the free Gamma case corresponds to $b^2 = 4c$, the free Poisson case to $c=0$, and the free Gaussian (semicircular) case to $b=c=0$. See also \cite{Boz-Bryc} for related results.

\br
Let $\phi$ be the free product state of $\set{\phi^{b_i, c_i}, i = 1, \ldots, n}$. The free cumulant generating function of $\phi$ is simply
\[
R(\mb{z}) = \sum_{i=1}^n R_i(z_i),
\]
where $R_i$ is the free cumulant generating function of $\phi^{b_i, c_i}$, satisfying
\[
R_i(z_i)/z_i^2 = 1 + b_i R_i(z_i)/z_i + c_i (R_i(z_i)/z_i)^2.
\]
Let $U_i(\mb{z}) = (R_i(z_i)/z_i)^{<-1>}$. Then the free Sheffer polynomials corresponding to $(R, \mb{U})$ are orthogonal. Indeed, these polynomials satisfy the recursion
\[
y_i P_{(j, \vec{u})}= P_{(i, j, \vec{u})},
\]
\[
y_i P_{(i, j, \vec{u})}= P_{(i, i, j, \vec{u)}} + b_i P_{(i, j, \vec{u})} + P_{(j, \vec{u})},
\]
\[
y_i P_{(i, i, \vec{u})} = P_{(i, i, i, \vec{u})} + b_i P_{(i, i, \vec{u})} + (1 + c_i) P_{(i, \vec{u})}.
\]
for $i \neq j$. So
\[
B_{i, (\alpha, \beta, \vec{w}), (s, t, \vec{u})} = \delta_{\vec{w}, \vec{u}} \delta_{si} \delta_{\alpha s} \delta_{\beta t} b_i
\]
and
\[
C_{i, (\alpha, \vec{w}), (s, t, \vec{u})} = \delta_{\vec{w}, \vec{u}} \delta_{si} \delta_{\alpha t} (1 + \delta_{ti} c_i).
\]
The conditions of Proposition~\ref{Prop:Orthogonal-recursion} are satisfied, so the polynomials are orthogonal.

\br
Explicitly, these polynomials are free products. Denote by $\set{P_k^{b, c}}$ the one-variable free Meixner polynomials from  Proposition~\ref{Prop:One-variable-free-Meixner}. Decompose a multi-index $\vec{u}$ so that
\[
x_{\vec{u}} = x_{v(1)}^{i(1)} x_{v(2)}^{i(2)} \ldots x_{v(k)}^{i(k)},
\]
where the consecutive indices $v(j) \neq v(j+1)$, although non-consecutive indices may coincide. Then
\[
P_{\vec{u}}(\mb{x}) = \prod_{j=1}^k P_{i(j)}^{b_{v(j)}, c_{v(j)}}(x_{v(j)}).
\]
\end{Ex}

\noindent
Thus free products of one-dimensional free Meixner states are free Meixner. The following proposition provides a partial converse.

\begin{Prop}
\label{Prop:Characterization}
Suppose that $\phi$ is a tracial free Meixner state with
\[
D_i D_j R_\phi = \delta_{ij} + \sum_t B_{ij}^t D_t R_\phi.
\]
Then up to a rotation, $\phi$ is a free product state of semicircular and free Poisson distributions.
\end{Prop}

\begin{proof}
It follows from Corollary~\ref{Cor:Commute} that the matrices $\set{B^r}$ are all symmetric and mutually commuting. So we can find an orthogonal transformation $O$ such that $(O^{-1} B^r O)_{ij} = \delta_{ij} b_i^r$ for all $r$. Performing the change of variable in Lemma~\ref{Lemma:Change}, we get
\[
O (Q R_{\mb{y}}(\mb{w})) O^{-1} = I + B \cdot O \mb{D} R_{\mb{y}}(\mb{w}).
\]
So
\[
Q R_{\mb{y}}(\mb{w}) = I + O^{-1} (B \cdot O \mb{D} R_{\mb{y}}(\mb{w})) O.
\]
Note that $B_{ij}^k = \sum O_{is} b_s^k O_{js}$ is also equal to $B_{ik}^j = \sum O_{is} b_s^j O_{ks}$. Then
\[
(O^{-1} (B \cdot O \mb{w}) O)_{\alpha \beta} = \sum_{k,l} \delta_{\alpha \beta} b_{\alpha}^k O_{kl} w_l.
\]
On the other hand, it is also equal to
\[
\sum_{i,j,s,k,l} O_{i \alpha} O_{is} b_s^j O_{ks} O_{j \beta} O_{kl} w_l = \sum_j b_{\alpha}^j O_{j \beta} w_{\alpha}.
\]
As a result, $\sum_k \delta_{\alpha \beta} b_\alpha^k O_{kl} = \sum_j \delta_{\alpha l} b_\alpha^j O_{j \beta}$ and
\[
(O^{-1} (B \cdot O \mb{w}) O)_{\alpha \beta} = \delta_{\alpha \beta} \bigl( \sum b_\alpha^k O_{k \alpha} \bigr) w_\alpha.
\]
Denote $b_\alpha = \sum_k b_\alpha^k O_{k \alpha}$. Then
\[
D_i D_j R_{\mb{y}}(\mb{w}) = \delta_{ij} + \delta_{ij} b_i D_i R_{\mb{y}}(\mb{w}) = \delta_{ij} \Bigl(1 + b_i D_i R_{\mb{y}}(\mb{w})) \Bigr).
\]
Therefore
\[
R_{\mb{y}}(\mb{w}) = \sum_{i=1}^n R_{y_i}(w_i),
\]
so all the mixed cumulants are zero and the components are freely independent. Moreover, each $R_i$ satisfies the equation
\[
R_i/w_i^2 = 1 + b_i R_i/w_i.
\]
This is exactly the equation in Proposition~\ref{Prop:One-variable-free-Meixner} for the free Poisson case, or for the semicircular case if $b_i = 0$.
\end{proof}

\section{A freely infinitely divisible example}
\label{Section:Example}
\begin{Defn}
A state $\phi$ is freely infinitely divisible if for all $t>0$, the functional $\phi^t$ with the free cumulant generating function
\[
R_{\phi^t}(\mb{z}) = t R_{\phi}(\mb{z})
\]
is also positive definite.
\end{Defn}

\begin{Remark}
One-dimensional free Meixner states $\phi^{b,c}$ of Proposition~\ref{Prop:One-variable-free-Meixner} and Example~\ref{Example:Product-states} are freely infinitely divisible for $c \geq 0$, and are not infinitely divisible for $-1 \leq c < 0$. In fact, in this case $\phi^t$ is a state only for $t \geq -(1/c)$.

\br
Thus all the states of Proposition~\ref{Prop:Characterization} are freely infinitely divisible, but some more general free product states of Example~\ref{Example:Product-states} are not. In this section, we construct an example of a freely infinitely divisible free Meixner state that is not a free product state.
\end{Remark}

\begin{Defn}
A functional $\psi$ on $\mf{R} \langle \mb{x} \rangle$ is conditionally positive definite if it is positive definite on polynomials of degree at least 2.
\end{Defn}

\begin{Lemma}
\label{Lemma:Generator}
$\phi$ is freely infinitely divisible if and only if its free cumulant functional is conditionally positive definite.
\end{Lemma}

\begin{proof}
$R_\phi[x_{\vec{u}}] = \frac{d}{dt}\bigl|_{t=0} M_{\phi^t}[x_{\vec{u}}]$. So if each $\phi^t$ is conditionally positive definite, so is $R_\phi$. For the converse, starting with a conditionally positive linear functional, one constructs symmetric operators with the joint distribution $\phi$. See \cite{GloSchurSpe} or Section 4 of \cite{AnsQCum}.
\end{proof}

\noindent
The following lemma is reminiscent of the Kolmogorov representation for infinitely divisible measures with finite variance.

\begin{Lemma}
\label{Lemma:CPD}
Let $\set{\phi_i, i = 1, \ldots, n}$ be positive definite functionals on $\mf{R} \langle \mb{x} \rangle$. Define the functional $\psi$ on $\mf{R} \langle \mb{x} \rangle$ as follows:
\begin{gather*}
\psi[1] = \psi[x_i] = 0, \\
\psi[x_i x_j] = \delta_{ij}, \\
\psi[x_i P(\mb{x}) x_j] = \delta_{ij} \phi_i[P(\mb{x})].
\end{gather*}
Then $\psi$ is conditionally positive definite.
\end{Lemma}

\begin{proof}
For such $\psi$,
\[
\begin{split}
\psi \Bigl[ \bigl( \sum_i P_i(\mb{x}) x_i \bigr)^\ast \bigl( \sum_j P_j(\mb{x}) x_j \bigr) \Bigr]
&= \psi \Bigl[ \bigl( \sum_i x_i P_i(\mb{x})^\ast \bigr) \bigl( \sum_j P_j(\mb{x}) x_j \bigr) \Bigr] \\
&= \psi \Bigl[ \sum_i x_i P_i(\mb{x})^\ast P_i(\mb{x}) x_i \Bigr] \\
&= \sum_i \phi_i[P_i(\mb{x})^\ast P_i(\mb{x})] \geq 0,
\end{split}
\]
so $\psi$ is conditionally positive definite.
\end{proof}

\noindent
We will denote $\psi$ as above by $\exp(\phi_1 \oplus \ldots \oplus \phi_n)$.

\br
The following result was already used in the proof of Theorem 3.21 of \cite{AnsAppell}; here we formulate it as a lemma. Considering how different the relation \eqref{Cumulant-moment} is from the logarithmic relation between moments and the usual cumulants, this result is surprisingly similar to the identity $(\log f)' = f'/f$.

\begin{Lemma}
\label{Lemma:Cumulant-moment}
For $z_i = w_i (1 + M(\mb{w}))$, we have
\[
(1 + M(\mb{w})) D_{z_i} R(\mb{z}) = D_{w_i} M(\mb{w}).
\]
\end{Lemma}

\begin{proof}
The result follows immediately from the relation \eqref{Cumulant-moment}.
\end{proof}

\begin{Prop}
Let $\psi$ be the distribution of a free semicircular system with means $b_i$ and variances $c_i$. In other words, $\psi$ is the state with the free cumulants
\[
R_\psi[x_i] = b_i, \qquad R_\psi[x_i^2] = c_i,
\]
and all the other free cumulants are zero. Define the state $\phi$ by $R_\phi = \exp(\psi^{\oplus n})$. Then $\phi$ is a free Meixner state.
\end{Prop}

\begin{proof}
By definition,
\begin{equation*}
R_{\psi}(\mb{z}) = \sum_i \Bigl( b_i z_i + c_i z_i^2 \Bigr).
\end{equation*}
So
\[
D_i R_{\psi}(\mb{z}) = b_i + c_i z_i.
\]
Using the change of variables $z_k = w_k \bigl(1 + M_{\psi}(\mb{w}) \bigr)$ and Lemma~\ref{Lemma:Cumulant-moment}, we get
\[
\bigl(1 + M_\psi(\mb{w}) \bigr)^{-1} D_i M_\psi(\mb{w}) = b_i + c_i w_i \bigl(1 + M_\psi(\mb{w}) \bigr),
\]
and so
\[
D_i M_\psi(\mb{w}) = b_i \bigl(1 + M_\psi(\mb{w}) \bigr) + c_i \bigl(1 + M_\psi(\mb{w}) \bigr) w_i \bigl(1 + M_\psi(\mb{w}) \bigr).
\]
Combination of Lemmas~\ref{Lemma:Generator} and~\ref{Lemma:CPD} shows that $\phi$ is a well-defined freely infinitely divisible state. Its free cumulant generating function is
\[
R_{\phi}(\mb{w}) = \sum_j w_j (1 + M_{\psi}(\mb{w})) w_j.
\]
Then
\[
D_j R_{\phi}(\mb{w}) = (1 + M_{\psi}(\mb{w})) w_j
\]
and
\[
\begin{split}
D_i D_j R_{\phi}(\mb{w})
&= \delta_{ij} + D_i M_{\psi}(\mb{w}) w_j \\
&= \delta_{ij} + b_i \bigl(1 + M_\psi(\mb{w}) \bigr) w_j + c_i \bigl(1 + M_\psi(\mb{w}) \bigr) w_i \bigl(1 + M_\psi(\mb{w}) \bigr) w_j \\
&= \delta_{ij} + b_i D_j R_{\phi}(\mb{w}) + c_i D_i R_{\phi}(\mb{w}) \ D_j R_{\phi}(\mb{w}).
\end{split}
\]
Thus $B_{ij}^t = \delta_{jt} b_i$, the conditions of Corollary~\ref{Cor:Orthogonal} are satisfied, and the free Sheffer polynomials corresponding to $\phi$ are orthogonal.
\end{proof}

\noindent
Note  that unless all $c_i = 0$, $\phi$ is not a tracial state.

\end{document}